\documentclass[12pt]{amsart}
\usepackage{amsmath,amssymb,amsthm}
\usepackage[]{latexsym}
\newtheorem{defn0}{Definition}[section]
\newtheorem{prop0}[defn0]{Proposition}
\newtheorem{thm0}[defn0]{Theorem}
\newtheorem{lemma0}[defn0]{Lemma}
\newtheorem{corollary0}[defn0]{Corollary}
\newtheorem{example0}[defn0]{Example}
\newtheorem{remark0}[defn0]{Remark}
\newtheorem{conjecture0}[defn0]{Conjecture}
\newtheorem{notation0}[defn0]{Notation}

\newenvironment{definition}{\begin{defn0}\rm}{\end{defn0}}
\newenvironment{proposition}{\begin{prop0}}{\end{prop0}}
\newenvironment{theorem}{\begin{thm0}}{\end{thm0}}
\newenvironment{lemma}{\begin{lemma0}}{\end{lemma0}}
\newenvironment{corollary}{\begin{corollary0}}{\end{corollary0}}

\newenvironment{remark}{\begin{remark0}\rm}{\end{remark0}}
\newenvironment{conjecture}{\begin{conjecture0}}{\end{conjecture0}}

\newcommand{\CM}{{\mathrm {CM}}}

\newcommand{\cT}{{\mathcal T}}

\newcommand{\om}{{\omega}}

\newcommand{\Gal}{{\mathrm {Gal}}}

\newcommand{\disc}{{\mathrm {disc }}}

\newcommand{\lra}{\longrightarrow}

\newcommand{\M}{\mathrm{M}}
\newcommand{\Aut}{\mathrm{Aut}}
\newcommand{\PGL}{{\mathrm{PGL}}}
\newcommand{\GL}{{\mathrm{GL}}}
\newcommand{\End}{{\mathrm{End}}}
\newcommand{\PSL}{{\mathrm {PSL}}}
\newcommand{\SL}{{\mathrm {SL}}}

\newcommand{\Z}{{\mathbb Z}}
\newcommand{\Q}{{\mathbb Q}}

\newcommand{\R}{{\mathbb R}}
\newcommand{\F}{{\mathbb F}}

\newcommand{\cO}{{\mathcal O}}
\newcommand{\tcO}{{\mathcal O^{(p)}}}

\newcommand{\n}{{\mathrm{n}}}

\newcommand{\ra}{{\rightarrow}}

\include{thebibliography}

\begin{document}

\title[Automorphisms on Shimura curves]{On the non-existence of exceptional automorphisms on Shimura curves}

\author{Aristides Kontogeorgis, Victor Rotger}

\address{Department of Mathematics, University of the \AE gean, 83200 Karlovassi, Samos, 
Greece\\ { \texttt{\upshape http://eloris.samos.aegean.gr}}}
\email{kontogar@aegean.gr}

\address{Escola Universit\`aria Polit\`ecnica de Vilanova i la
Geltr\'u, Av. V\'{\i}ctor Balaguer s/n, E-08800 Vilanova i la Geltr\'u, Spain.}
\email{vrotger@ma4.upc.edu}



\begin{abstract}
We study the group of automorphisms of Shimura curves $X_0(D, N)$ attached to an Eichler order of square-free level $N$ in an indefinite rational quaternion algebra of discriminant $D>1$. We prove that, when the genus $g$ of the curve is greater than or equal to $2$, $\Aut (X_0(D, N))$ is a $2$-elementary abelian group which contains the group of Atkin-Lehner involutions $W_0(D, N)$ as a subgroup of index $1$ or $2$. It is conjectured that $\Aut (X_0(D, N )) = W_0(D, N)$ except for finitely many values of $(D, N)$ and we provide criteria that allow us to show that this is indeed often the case. Our methods are based on the theory of complex multiplication of Shimura curves and the Cerednik-Drinfeld theory on their rigid analytic uniformization at primes $p\mid D$.
\end{abstract}

\maketitle

\section{The automorphism group of Shimura curves}

\subsection{Congruence subgroups of $\PSL_2(\R )$ and automorphisms}

$\\ $

Let $\Gamma $ be a congruence subgroup of $\PSL_2(\R )$. As explained in \cite[\S 4]{LMR}, $\Gamma $ is a congruence subgroup of $\PSL_2(\R )$ if there exists

\begin{itemize}
\item A quaternion algebra $B/F$ over a totally real number field $F$ of degree $d\geq 1$,
\item An embedding $\varphi: B\,\hookrightarrow \M_2(\R )\times D\times \stackrel{(d-1)}{...} \times D$,
\item An integral two-sided ideal $I$ of a maximal order $\cO $ of $B$,
\end{itemize} 
such that $\Gamma $ contains $\varphi (\{ \alpha \in \cO^1: \alpha \in 1+I\})$.

\vspace{0.4cm}

Here, we let $D$ denote Hamilton's skew-field over $\R $ and $\n: B \lra F$ stand for the reduced norm. We write $\cO^1=\{ \alpha \in \cO: \n(\alpha )=1\}$. We refer the reader to \cite{Vi} for generalities on quaternion algebras. Examples of congruence subgroups of $\PSL_2(\R )$ with $F=\Q$ will be described in detail below.

Let $X_{\Gamma }$ denote the compactification of the Riemann surface $\Gamma \backslash \mathcal H$. Let $N = \mathrm{Norm}_{\PSL_2(\R )}(\Gamma)$ denote the  normalizer of $\Gamma $ in $\PSL_2(\R )$. The group $B_{\Gamma } = N /\Gamma $ is a finite subgroup of $\Aut(X_{\Gamma })$. 

When the genus of $X_{\Gamma }$ is $0$ or $1$, $\Aut(X_{\Gamma })$ is not a finite group and necessarily $\Aut(X_{\Gamma }) \varsupsetneq B_{\Gamma }$. However, there exist finitely many congruence groups $\Gamma $ for which $g(X_{\Gamma })\le 1$. One actually expects much more, as we claim in the following conjecture. 

\begin{conjecture}\label{Conj} $\Aut(X_{\Gamma }) = B_{\Gamma }$ for all but finitely many congruence groups $\Gamma \subset \PSL_2(\R )$. 
\end{conjecture}

We call {\em exceptional} those congruence groups $\Gamma $ for which the genus of $X_{\Gamma }$ satisfies $g\geq 2$ and $\Aut (X_{\Gamma })\varsupsetneq B_{\Gamma }$. 

The conjecture as we have stated remains widely open, although it is motivated by several positive partial results towards it. The first general statement is the following classical result of Riemann surfaces, of which we briefly recall a proof.

\begin{proposition}\label{cover} Let $\Gamma $ be a congruence subgroup of $\PSL_2(\R )$ that contains no elliptic nor parabolic elements. Then $\Aut(X_{\Gamma }) = B_{\Gamma }$.
\end{proposition}

{\bf Proof.} Since $\Gamma $ contains no parabolic elements, the quotient $\Gamma \backslash \mathcal H$ is already compact. The absence of elliptic elements implies that the natural projection $\mathcal H\, \lra \, X_{\Gamma }=\Gamma \backslash \mathcal H$ is the universal cover of the curve. Thus all automorphisms of $X_{\Gamma }$ lift to a M\"{o}bius transformation of $\mathcal H$ which by construction normalizes $\Gamma $. The result follows. $\Box $

Besides this, the question has been settled for certain families of modular curves, as we now review.

Let $D\ge 1$ be the square-free product of an even number of
prime numbers and let $N\ge 1$, $(D, N)=1$, be an integer coprime to $D$. 

Let $B$ be a quaternion algebra over $\Q $ of reduced discriminant $D$ such that there exists an monomorphism $B \, \stackrel{\varphi }\longrightarrow  \, \mathrm{M}_2(\R )$. Let $\cO $ be a maximal order in $B$. Regard $\cO^1$ as a subgroup of $\SL_2(\R )$ by means of $\varphi $.

For any prime $p\mid N$ fix isomorphisms $\cO \otimes \Z_p \,\simeq \,\mathrm{M}_2(\Z_p)$. If $p^{e}\mid \mid N$ is the exact power of $p$ which divides $N$, let $\pi_p: \, \cO\otimes \Z_p\,\simeq \, \mathrm{M}_2(\Z_p)\, \rightarrow \, \mathrm{M}_2(\Z /p^{e}\Z )$ denote the natural reduction map.

Define the congruence groups $\Gamma(D, N) \subseteq \Gamma_1(D, N) \subseteq \Gamma_0(D, N)$ as follows:

\vspace{0.3cm}

\begin{itemize}

\item $\Gamma(D, N) = \{ \gamma \in \cO^1: \, \pi_p (\gamma ) = \mathrm{Id} \in \GL_2(\Z /p^{e}\Z )\mbox{ for all }p\mid N\}$

\vspace{0.3cm}

\item $\Gamma_1(D, N) = \{ \gamma \in \cO^1: \, \pi_p (\gamma ) = \begin{pmatrix}
1& * \\ 0 & 1 
\end{pmatrix}\in \GL_2(\Z /p^{e}\Z )\mbox{ for all }p\mid N\} $

\vspace{0.3cm}

\item $\Gamma_0(D, N) = \{ \gamma \in \cO^1: \, \pi_p (\gamma ) = \begin{pmatrix}
*& * \\ 0 & * 
\end{pmatrix}\in \GL_2(\Z /p^{e}\Z )\mbox{ for all }p\mid N\}.
$

\end{itemize}

Let $X(D, N)$, $X_1(D, N)$ and $X_0(D,N)$ denote the corresponding Shimura modular
curve of discriminant $D$ and level $N$ (cf.\,\cite{BoCa}). 

When $D=1$, these are just another names for the elliptic modular curves $X(N)$, $X_1(N)$ and $X_0(N)$ which classify elliptic curves with various level structures. When $D>1$, these curves still admit a moduli interpretation in terms of abelian surfaces with quaternionic multiplication (cf.\,\cite{BoCa}). 

Finally, note that when $N=1$ we have $X(D, 1) = X_1(D, 1) = X_0(D, 1)$ and we shall simply denote this curve by $X_D$. 

As further examples of congruence subgroups of $\PSL_2(\R )$ are the normalizers $B_{\Gamma }$ for the groups $\Gamma $ above.

The group $B_{\Gamma_0(D, N)}$ contains the subgroup $W_0(D, N)$ of Atkin-Lehner modular involutions (cf.\,\cite{Ogg:77}, \cite{Ogg1}). It can be described as
$$W_0(D, N) = \{ \om_m : m\ge 1, m\mid D\cdot N, (m, D N/m) = 1\},$$ where $\om_m^2 = w_1 = \mathrm{Id}$ and $\om_m\cdot \om_{m'} = \om_{m m'/(m, m')^2}$. Hence $$W_0(D, N) \simeq (\Z /2\Z )^{|\, \{ p\mid D\cdot N \} \,|}.$$ 

If $N$ is square-free we actually have $B_{\Gamma_0(D, N)} = W_0(D, N)$. 

\begin{remark}\cite[p.\,1]{LMR}
The above conjecture cannot be true for the wider family of arithmetic subgroups of $\PSL_2(\R )$. As an example for which the conjecture fails, let $X(2)$ be the classical modular curve of full level structure $\Gamma(2)$. This is a rational curve with three cusps and no elliptic points. If we choose a rational coordinate $z$ for which the cusps are at $0$, $1$ and $\infty$, then $z^{1/n}$ is a rational coordinate for a subgroup of index $n$ in X(2). This subgroup is not congruence for $n$ large enough, but it is arithmetic. The corresponding curve has genus $0$ and the automorphism group is not finite.
\end{remark}

Several results towards Conjecture \ref{Conj} regarding the above families of curves have been achieved by different authors. 

\begin{theorem}\label{all} For the curves below, if their genus is greater than or equal to two, one has:

\begin{enumerate}

\item[(i)] (Serre \cite{Se}, Kontogeorgis \cite{KonPHD}) $\Aut(X(N)) = B_{\Gamma(1, N)}$ for all $N$.

\vspace{0.3cm}

\item[(ii)] (Momose \cite{Mo}, Buzzard \cite{Bu}) $\Aut(X_1(D, N)) = B_{\Gamma_1(D, N)}$ for all $N\geq 4$.

\vspace{0.3cm}

\item[(iii)] (Ogg \cite{Ogg:77}, Kenku-Momose \cite{KeMo}, Elkies \cite{El}), $\Aut(X_0(N)) = B_{\Gamma_0(1, N)}$ for all $N\ne 37, 63$.

\vspace{0.3cm}

\item[(iv)] (Baker-Hasegawa \cite{BaHa}) $\Aut(X_0(p)/\langle \om_p \rangle )$ is trivial for all primes $p\ne 67, 73, 103,\\ 107, 167, 191$. 

\vspace{0.3cm}

\item[(v)] (Rotger \cite{Rot:2002}) $\Aut (X_D) = W_0(D, 1)$, where $D=2 p$ or $3 p$ for some prime $p$.

\end{enumerate}

\end{theorem}

Part $(i)$ was shown by Serre \cite{Se} for prime $N$ and extended to composite $N$ by Kontogeorgis in \cite{KonPHD}.

Part $(ii)$ above has been proved by Momose \cite{Mo} when $D=1$. When $D>1$, $\Gamma_1(D, N)$ contains no parabolic (because $B$ is division) nor elliptic elements by \cite[Lemma 2.2]{Bu}. Thus Proposition \ref{cover} applies.

The analogous result of $(i)$ in positive characteristic for the Drinfeld modular curves $X(\mathfrak N)$ has been achieved in \cite{CKK}. 

In $(iv)$, for $p=67, 73, 103, 107, 167, 191$ the genus of $X_0(p)/\langle \om_p \rangle $ is $2$ and the hyperelliptic involution is the only (exceptional) automorphism of the curve.

\subsection{Main results}

$\\ $

Let $D>1$ be the square-free product of an even number of primes and $N\ge 1$ be a square-free integer coprime to $D$. Let $r= | \{ p\mid D N\} |$. 

The aim of this section is making progress towards Conjecture \ref{Conj} for the family of Shimura curves $X_0(D, N)$. Due to their moduli interpretation, these curves admit a canonical model over $\Q $ which we still shall denote $X_0(D, N) /\Q  $. There exists a flat proper model $M_0(D, N)/\Z $ of $X_0(D, N)$ which extends the moduli interpretation to arbitrary schemes over $\Z $ and is smooth over $\Z [\frac{1}{D N}]$. 

For primes $p\nmid D$, the construction of this model over $\Z_p$ is very similar to that of the elliptic modular curve $X_0(N)$ as in \cite{DeRa}; see \cite{Bu} for more details. For primes $p\mid D$ the description of $M_0(D, N)\otimes \Z_p$ is due to Cerednik and Drinfeld. For any prime $p$ we let $M_0(D, N)_p$ denote the closed fiber of $M_0(D, N)$ at $p$.

\begin{proposition}\label{comp} Let $U \subseteq W_0(D, N)$ be a subgroup and let $X_0(D, N)/U$ denote the quotient curve. If the genus of $X_0(D, N)/U$ is at least $2$, then all automorphisms of $X_0(D, N)/U$ are defined over $\Q $ and $$\Aut (X_0(D, N)/U) \simeq (\Z/2\Z )^s$$ for some $s\ge r - \mathrm{rank}_{\F_2}(U)$.
\end{proposition}

{\bf Proof.} Write $X=X_0(D, N)$. Let $J_U$ and $J$ denote the Jacobian varieties of $X/U$ and $X$, respectively. A result essentially due to Ribet claims that $J$ is isogenous over $\Q $ to $J_0(D\cdot N)^{D-new}$, the $D$-{\em new part} of the Jacobian $J_0(D, N )$ of $X_0(D N )$; see e.\,g.\,\cite[Theorem 5.4]{GoRo2} or \cite{Ri} for more details. Since $D N$ is square-free, it is well-known that $J_0(D N)$ has semi-stable reduction over $\Q $ (cf.\,\cite{DeRa}) and $\End_{\Q }(J_0(D N))\otimes \Q $ is a product of totally real fields (cf.\,\cite{Ri75}). The same thus holds for $J$ and $J_U$. The proposition now follows as a direct application of \cite[Proposition 2.4]{BaHa}. $\Box $

\vspace{0.3cm}

Let us fix some notation. Throughout, let $s\ge r$ be the integer such that $|\Aut(X_0(D, N))| = 2^s$. We have $s\geq r$ and Conjecture \ref{Conj} in this setting predicts that $s=r$. 

Letters $p, q$ will stand for non-necessarily different prime numbers and $(\frac{\,\cdot \,}{p})$ shall denote the Kronecker quadratic character {\em mod} $p$.
For an imaginary quadratic field $K$, let $\delta _K$ denote its discriminant and $h(K)$ its class number. 
If $\delta_R=\disc(R)$ is the discriminant of some order $R$ in $K$, let $h(\delta_R)$ denote its class number. 

If $m>1$ is a square-free integer, let $\delta_m = -4$ if $m=2$; $\delta_m = \delta_{\Q(\sqrt{-m})}$ otherwise.  

For any prime $p\mid D N$, set $\varepsilon_p = \begin{cases}
1 & \text{ if } p\mid D \\
-1 & \text{ if } p \mid N
\end{cases}$ 

and for any $m\mid D N$, 
$$
h_{D, N}(m) = \begin{cases}
1 & \text{ if } m=2 \\
h(-4 m) & \text{ if } m\ne 2, m\not \equiv 3 \text{ mod } 4 \\
h(-m) & \text{ if } m\equiv 3 \text{ mod } 4 \text{ and } h(-4 m)>h(-m) \text{ or } 2\mid D N \\
2 h(-m) & \text{if } m\equiv 3 \text{ mod } 4, h(-4m )=h(-m), 2\nmid D N.
\end{cases}$$

If in addition $(\frac{\delta_m}{p})\ne \varepsilon_p$ for all $p\mid D N$ if $m\ne 2$; $(\frac{-4}{p})\ne \varepsilon_p$ for all $p\mid D N$ or $(\frac{-2}{p})\ne \varepsilon_p$ for all $p\mid D N$ if $m=2$, set
$$
\sigma_{D, N}(m)=\begin{cases}
|\{ p\mid D N, (\frac{\delta_m}{p}) = -\varepsilon_p\} |  & \text{ if } m\ne 2,\\
\text{min }(|\{ p\mid D N, (\frac{-4}{p}) = -\varepsilon_p\} |, |\{ p\mid D N, (\frac{-2}{p}) = -\varepsilon_p\} |) & \text{ if } m=2.    \\
\end{cases}
$$

Otherwise, set $\sigma_{D, N}(m)=\infty $.

\begin{theorem}\label{CM} Assume $g(X_0(D, N))\geq 2$.
\begin{enumerate}
\item[(i)] If there exist $p, q\mid D N$ such that $(\frac{-4}{p}) = \varepsilon_p$ and $(\frac{-3}{q}) = \varepsilon_q$, then $s=r$.
\item[(ii)] Let $m\mid D N$ such that $\sigma_{D, N}(m)<\infty $. Then $$s\leq \text{ord}_2\,( h_{D, N}(m) )+ \sigma_{D, N}(m) +1.$$ 
\item[(iii)] $s\leq \text{ord}_2 (g-1) + 2$. 
\end{enumerate}
\end{theorem}

For any two coprime square-free integers $\delta , \nu \ge 1$, let $h(\delta , \nu )$ denote the class number of any Eichler order of level $\nu $ in a quaternion algebra of reduced discriminant $\delta $ over $\Q $, which counts the number of one-sided ideals of the order up to principal ideals. An explicit formula for $h(\delta , \nu )$ is given in \cite[p.\,152]{Vi}.

\begin{theorem}\label{CD} Assume $g(X_0(D, N))\geq 2$.
\begin{enumerate}
\item[(i)] Let $m=2$ or $3$. Assume $(\frac{\delta_m}{p}) \ne \varepsilon_p$ for all $p\mid D N$ except at most for one prime divisor of $D$. If $m\mid D N$ then $s=r$; otherwise $s\leq r+1$.

\item[(ii)] For any odd $p\mid D$: $s\leq \text{ord}_2 \, h(D/p, N) + 3$. 

\item[(iii)] For any $\ell \nmid 2 D N$: $s\leq \text{ord}_2 \, |\,M_0(D, N)_{\ell }(\F_{\ell })\,|\,+\,1$.
\end{enumerate}
\end{theorem}

An explicit formula for the number of rational points of $M_0(D, N)_{\ell }$ over $\F_{\ell ^n}$, $n\geq 1$, may be found in \cite[\S 2]{RSY}. As a direct consequence of Theorem \ref{CM} $(i)$ and Theorem \ref{CD} $(i)$ we derive the following result.

\begin{corollary}\label{Shimurabound} If $g(X_0(D, N))\geq 2$ then $s= r$ or $r+1$.
\end{corollary}

We prove Theorem \ref{CM} in Section \ref{CCMM}, while the proof of Theorem \ref{CD} is offered in Section \ref{CCDD}.

Theorems \ref{CM} and \ref{CD} may be applied to show that there exist no automorphisms on many Shimura curves beyond the Atkin-Lehner involutions. Indeed, Theorem \ref{CM} $(ii)$ covers the case $e_2=e_3=0$ while Theorem \ref{CD} $(ii)$ solves the cases $e_2\ne 0$, $2\mid D N$ and $e_3\ne 0$, $3\mid D N$. As a consequence, we have the following corollary.

\begin{corollary}\label{corol}
If $6\mid D N$ then $s=r$, that is, $\Aut(X_0(D, N)) = W_0(D, N)$.
\end{corollary}

All together, Theorems \ref{CM} and \ref{CD} can be applied to many other Shimura curves. It follows for instance from Theorem \ref{CM} $(iii)$ that $\Aut (X_0(2 p, N))=W_0(2 p, N)$ for all primes $p$ and $N$ are primes, $N\equiv 3$ mod $8$. This follows because $h(2, N)$ is always odd (cf.\,\cite[p.\,152]{Vi}). We refer the reader to Proposition \ref{2000} for more numerical computations.

\subsection{Overview of the article}

We devote the remainder of this note to introduce the necessary tools which we shall need to eventually prove Theorems \ref{CM} and \ref{CD}.

Next section reviews the theory of complex multiplication on Shimura curves and its behaviour with respect to the Atkin-Lehner group. Following ideas borrowed from \cite{Rot:2002}, we use this to offer a proof of Theorem \ref{CM}.

Section \ref{CCDD} recalls the theory of Cerednik-Drinfeld, which provides an explicit description of the reduction mod $p$ of Shimura curves $X_0(D, N)$ at primes $p\mid D$. 

Namely, the main result of Cerednik and Drinfeld describe $X_0(D, N)\times \Q_p$ as a quadratic twist of a Mumford curve over $\Q_p$ which is rigid analytically uniformized by a certain discrete finitely generated subgroup $\Gamma_+$ of $\PGL_2(\Q_p)$. The group $\Gamma_+$ is constructed by means of an {\em interchange of invariants} of the quaternion algebra of reduced discriminant $D$ over $\Q$.

In turn, this allows us to interpret the dual graph of the special fiber of a suitable model of $X_0(D, N)$ over $\Z_p$ as the quotient of the Bruhat-Tits tree at $p$ by $\Gamma_+$.

We use this to give a combinatorial description of the stable model and minimal regular model of $X_0(D, N)$ at primes $p\mid D$: see Proposition \ref{me}.

In subsection \ref{pcd} we make use of this material to complete the proof of Theorem \ref{CD}. Finally, in subsection \ref{nr} we offer a numerical result which shows how the methods of this note apply to prove the non-existence of exceptional automorphisms on most Shimura curves $X_0(D, N)$ for $D\leq 1500$ and $N=1$.

\section{Automorphisms and points of complex multiplication}\label{CCMM}

Let $(D, N)$ be a pair as in the previous section. For an order $R$ in an imaginary quadratic field $K$, let $c_R$ be its conductor in $K$ and let $\mathrm{CM}(\delta_R)$ denote the set of complex multiplication (CM) points on $X_0(D, N)$ by the order $R$. A fundamental result of Shimura states that the coordinates of a CM-point $P\in \mathrm{CM}(\delta_R)$ on $X_0(D, N)$ generate the ring class field $H_R$ over $K$ (cf.\,e.\,g.\,\cite[\S 5]{GoRo2}). That is, 

\begin{equation}\label{Shi}
K\cdot \Q (P)=H_R.  
\end{equation}

The cardinality of $\mathrm{CM}(\delta_R)$ is given in \cite[Section 1]{Ogg1}:

\begin{equation}\label{CCM}
|\,\mathrm{CM}(\delta_R)\,|= \begin{cases}
0\quad  \text{if }(\frac{\delta_K}{p})= \varepsilon_p \text{ for some }p\mid D N \text{ or } (c_R, D N) \ne 1 \\
h(R)\cdot 2^{| \{ p\mid D N, (\frac{\delta_K}{p}) = -\varepsilon_p\} |}\quad  \text{otherwise.} \end{cases}
\end{equation}

CM-points arise in a natural way as fixed points of Atkin-Lehner involutions on $X_0(D, N)$. Indeed, for any $m\mid D N$, the set $\mathcal F_m$ of fixed points of $\om_m$ on $X_0(D, N)$ is 
\begin{equation}\label{FP}
\mathcal F_m = \begin{cases}
\mathrm{CM}(-4) \cup \mathrm{CM}(-8) & \text{if }m=2 \\
\mathrm{CM}(-m) \cup \mathrm{CM}(-4 m)& \text{if }m\equiv 3 \text{ mod }4\\
\mathrm{CM}(-4 m) & \text{otherwise.}
\end{cases}.
\end{equation}

\vspace{0.5cm}

Under our assumptions on $D$ and $N$, the groups $\Gamma_0(D, N)$ contain no parabolic elements and the only elliptic points on these curves are of order $2$ or $3$. In fact, the set of elliptic points of order $i=2, 3$ is $\CM (\delta_i)$. Thus, their cardinality $e_i$ is

\begin{equation}\label{ek}
e_i=\begin{cases}
0 & \text{if there exists } p\mid D N, (\frac{\delta_i}{p})=\varepsilon_p, \\
2^{r-1} & \text{if } i\mid D N \text{ and } (\frac{\delta_i}{p})\ne \varepsilon_p \text{ for any }p\mid D N, \\
2^{r} & \text{otherwise.}
\end{cases}
\end{equation}

By \cite[p.\,280,
301]{Ogg1} the genus of $X_0(D, N)$ is $$
g=g(X_0(D, N))=1+\frac{D N}{12}\cdot \prod_{p\mid D}(1-\frac{1}{p})\cdot
\prod_{p\mid N}(1+\frac{1}{p})-\frac{e_3}{3}-\frac{e_2}{4}.
$$

Next lemma is a particular case of \cite[\S 1, Hilfsatz 1]{Ogg:77}.

\begin{lemma}\label{Og}
Let $X$ be an irreducible curve over a field $k$ with $\mathrm{char} (k)\ne 2$. If $\Aut (X)\simeq (\Z /2 \Z )^s$ for some $s\geq 1$ and $P\in C(k)$ is a regular point, then $\mathrm{Stab }(P) = \{ \om \in \Aut (X) : \om (P)=P\}$ has order $1$ or $2$.
\end{lemma}

\vspace{0.3cm}

\subsection{Proof of Theorem \ref{CM}.}\label{pcm} Let $X=X_0(D, N)$ and $A=\Aut(X_0(D, N))$.

$(i)$ By (\ref{ek}) the group $\Gamma_0(D, N)\subset \PSL_2(\R )$ has no elliptic nor parabolic elements. According to Proposition \ref{cover}, $A = B_{\Gamma_0(D, N)} = W_0(D, N)$. 

$(ii)$ As it is clear from (\ref{CCM}) and (\ref{FP}), the assumptions of $(ii)$ imply that the set $\mathcal F_m$ of fixed points of $\om_m$ is non-empty. 

Since all automorphisms of $X$ commute with $\om_m$ by Proposition \ref{comp}, we deduce that $A$ acts on $\mathcal F_m$. 

Assume that either $m=2$ or $m\equiv 3$ mod $4$ and $h(-4 m)>h(-m)$ or $2\mid D N$. Set $S_1=\mathrm{CM}(-4)$, $S_2=\mathrm{CM}(-8)$ if $m=2$; $S_1=\mathrm{CM}(-m)$, $S_2=\mathrm{CM}(-4 m)$ otherwise. 

By (\ref{FP}), $\mathcal F_m = S_1\cup S_2$. Moreover, (\ref{CCM}) guarantees that at least one of $S_1$ and $S_2$ is non-empty.
In fact, when $m\equiv 3$ mod $4$, we have $S_1\ne \emptyset $. 

When $m=2$, $S_1\subseteq X(\Q(\sqrt{-1}))$ and $S_2\subseteq X(\Q(\sqrt{-2}))$ by (\ref{Shi}). As all automorphisms of $X$ are defined over $\Q $ by Proposition \ref{comp}, $A$ leaves both $S_1$ and $S_2$ invariant. 

When $m\equiv 3$ mod $4$, by (\ref{Shi}) any point in $S_1$ generates the Hilbert class field of $K=\Q(\sqrt{-m})$, which is an abelian extension of $K$ of degree $h(-m)$. Similarly, $S_2 = \emptyset $ if $2\mid D N$ by (\ref{CCM}); otherwise, any point in $S_2$ generates an extension of $K$ of degree $h(-4m)$. 

Since $h(-4 m)>h(-m)$, we conclude as above that $A$ fixes the sets $S_1$ and $S_2$.

Hence, in any case, $A$ acts on a non-empty set $S \,(=S_1$, $S_2$ or $\mathcal F_m)$ with $| S | = h_{D, N}(m)\cdot 2^{\sigma_{D, N}(m)}$. By Lemma \ref{Og} the stabilizer of any of the elements of $S$ in $A$ is exactly $\langle \om_m\rangle $. Thus $A/\langle \om_m\rangle $ acts freely on $S$ and $(ii)$ follows.
  
$(iii)$ Let $Y=X/A$ and $\pi : X \,\ra \, Y$ be the natural projection map. By Proposition \ref{comp}, $\pi $ is a finite morphism of degree $2^s$ which ramifies precisely at the set $\mathcal F$ of all fixed points of automorphisms of $X$. Riemann-Hurwitz's formula applied to $\pi $ says that $g(X) - 1 = 2^s (g(Y) - 1) + \frac{1}{2}\cdot | \mathcal F |$. Hence, it suffices to show that $\mathrm{ord}_2 (| \mathcal F |) \geq s-1$. But this fact readily follows from an induction argument, since for any $u \in A$ one has $| \mathcal F | = 2\cdot | \mathcal F' |$, where $\mathcal F' = \{ [x]\in X/\langle u \rangle : \om ([x]) = [x] \text{ for some } \om \in A/\langle u \rangle  \}$. $\Box $

\section{Cerednik-Drinfeld theory}\label{CCDD}

Fix a prime $p\mid D$. Let $k_p$ denote the quadratic unramified extension of $\Q_p$, $\Q_p^{unr}$ the maximal unramified extension of $\Q_p$ and $\Z_p^{unr}$ its ring of integers.

Let $\tcO $ be an Eichler order of level $N$ in a definite quaternion algebra of discriminant $D/p$ and fix an immersion $\tcO \hookrightarrow \mathrm{M}_2(\Q_p)$. Let $$\Gamma_+ = \{ \gamma \in (\tcO\otimes \Z[1/p])^* : \text{ord}_p (\det(\gamma ))\in 2 \Z \} / \Z[1/p]^* \hookrightarrow \PGL_2(\Q_p),$$ which is a finitely generated discontinuous subgroup of $\PGL_2(\Q_p)$.

We warn the reader that $\Gamma_+$ may not be torsion free: for $m=2, 3$, any embedding $\Z [\sqrt{\delta_m}] \stackrel{\varphi }{\hookrightarrow }\tcO $ produces an element $\gamma = \varphi (\sqrt{\delta_m})$ in $\Gamma_+$ of order $m$. However, by \cite[p.\,19]{GP} there exists a torsion-free normal subgroup $\Gamma_+^0 $ of finite index in $\Gamma_+$. The group $\Gamma_+^0$ is thus a Schottky group; let $X_{\Gamma^0_+} = \Gamma^0_+ \backslash (\mathbb{P}^1_{\Q_p} - \mathcal{L}_{\Gamma^0_+})$ denote the Mumford curve over $\Q_p$ attached to $\Gamma_+^0$ as in \cite{GP} or \cite{Mumford:72}. If its genus $g$ is at least $2$, $X_{\Gamma^0_+}$ has stable totally split reduction over $\Q_p$.

Since $A=\Gamma_+^0 \backslash \Gamma_+$ is a finite group which naturally lies in $\Aut (X_{\Gamma_+^0})$, there exists an algebraic curve $X_{\Gamma_+}$ over $\Q_p$ which is the quotient of $X_{\Gamma_+^0}$ by $A$.

\begin{theorem}[Cerednik-Drinfeld]

Let $\chi : \Gal(k_p/\Q_p) \,\rightarrow \Aut(X_{\Gamma_+}\otimes k_p)$, $\mathrm{Fr}\,\mapsto \, \om_p$ and $X_{\Gamma_+}^{\chi }$ be the quadratic twist of $X_{\Gamma_+}$ by $\chi $. Then

$$X_0(D, N)\times \Q_p\simeq X_{\Gamma_+}^{\chi }.$$
\end{theorem}

We refer the reader to \cite{BoCa} for a proof. See also \cite{JoLi}.

Let $\cT $ denote the Bruhat-Tits tree attached to $\Q_p$. Following \cite[\S 4, \S 5]{K} and \cite[\S 3]{JoLi}, the special fiber of $M_0(D, N)\otimes \Z_p$ is described up to a quadratic twist by the finite graph $\mathcal G=\Gamma_+ \backslash \mathcal T$, regarded as a {\em graph with lengths}. 

Each vertex $v$ and edge $e$ of $\mathcal G$ is decorated with the order $\ell(v)$ (resp.\,$\ell(e)$) of the stabilizer of $v$ (resp.\,$e$) in $\Gamma _+$, which we call its length. Geometrically, a vertex $v$ corresponds to an irreducible rational component $C_v$ of $M_0(D, N)_{p}$. An edge $e$ of length $\ell $ joining $v$ and $v'$ corresponds to an intersection point $P_e$ of $C_v\cap C_{v'}$ locally at which the scheme $M_0(D, N)\times \Z_p^{unr}$ is isomorphic to $\mathrm{Spec} (\Z_p^{unr}[X, Y]/(X Y -p^{\ell }))$. In particular, any automorphism of $X_0(D, N)$ induces an automorphism of $\mathcal G$ which leaves the set of edges of given length invariant.

Let $h=h(D/p, N)$. The number of vertices of $\mathcal G$ is\footnote{According to formula $(4.1)$ in \cite{K} and its notation, the number of vertices of $\Gamma_0\backslash \Delta $ is $h(D/p, 1)$. Since the index of $\Gamma_+$ in $\Gamma_0$ is $2$, the number of vertices of $\Gamma_+\backslash \Delta $ is $2 h(D/p, 1)$. These formulas extend to the case of non-trivial level $N$ without difficulty.} $2 h$ and that of edges is $h(D/p, N p)$. 
The set $\mathrm{Ver}(\mathcal G)$ of vertices of $\mathcal G$ may be written as 

\begin{equation}\label{vertices}
\mathrm{Ver}(\mathcal G) = V\cup V',\,\, V=\{ v_1, ..., v_h\},\, V'=\{ v_1', ..., v_h'\}, 
\end{equation} 
in such a way that the Atkin-Lehner involution $\om_p\in \Aut(X_0(D, N))$ acts on $\mathcal G$ as $\om_p(v_i) = v_i'$. There are no edges in $\mathcal G$ joining two vertices from the same set $V$ or $V'$ and hence  there are no loops in $\mathcal G$. We have $\ell(v_i) = \ell(v_i')$ and the number of edges of given length joining a vertex $v_i$ with $v_j'$ coincides with that of $v_i'$ with $v_j$. 

For a vertex $v$, it holds that $\ell(e)\mid \ell(v)$ for all edges $e$ in its star and 

\begin{equation}\label{p+1}
\sum_{e\in \mathrm{Star}(v)} \frac{\ell(v)}{\ell(e)} = p+1.
\end{equation}

When $(D/p, N) = (2, 1)$ or $(3,1)$, $\mathcal G$ consists of two vertices $v, v'$ of length $12$ ($6$, respectively), joined by $g+1$ edges. Assume otherwise that $(D/p, N)\ne (2, 1), (3, 1)$. Then all lengths of vertices and edges are $1$, $2$ or $3$. By \cite[(4.2)]{K}, for $\ell =2, 3$, the cardinality $h_{\ell }$ of vertices in $V$ of length $\ell $ is 
\begin{equation}\label{h}
h_{\ell } = \frac{1}{2} \prod_{q \mid \frac{D N}{p}}(1-\varepsilon_q\cdot (\frac{\delta_{\ell }}{q})).
\end{equation}

If $v$ is a vertex of length $\ell =2$ or $3$, by \cite[Proposition 4.2]{K} the number of edges of length $\ell $ in its star is 
\begin{equation}\label{star}
|\mathrm{Star}_{\ell }(v)| = 1+(\frac{\delta_{\ell }}{p}) \in \{ 0, 1, 2\}.
\end{equation}

The scheme $M_0(D, N)$ is regular if and only if $\ell (e) =1$ for all edges of $\mathcal G$. In general, a desingularization $\tilde{M}_0(D, N)$ of $M_0(D, N)$ is obtained by blowing-up $\ell (e) - 1$ times each singular point $P_e$. The resulting dual graph $\tilde{\mathcal G}$ is constructed from $\mathcal G$ by replacing each edge $e$ of length $\ell(e)\ge 2$, by a chain of $\ell (e)$ edges of length $1$ each (cf.\,\cite[Proposition 3.6]{JoLi}). 

In general, $M_0(D, N)$ is neither a minimal regular model nor a stable model of $X_0(D, N)$. Next proposition, which may be of independent interest,  shows how to construct these two models.

\begin{proposition}\label{me}
Assume $g=g(X_0(D, N))\geq 1$.

\begin{enumerate}

\item[(i)] Let $M_0(D, N)_{min }$ denote the minimal regular model of $X_0(D, N)$.
\begin{itemize}
\item If $p>2$ or $h_3=0$, $M_0(D, N)_{min } = \tilde{M}_0(D, N)$.

\item If $p=2$ and $h_3\geq 1$, $M_0(D, N)_{min }$ is the blow-down of all components $C_v$ of $\tilde{M}_0(D, N)$ with $\ell (v)=3$. 

Its dual graph $\mathcal G_{min}$ is obtained from $\tilde{\mathcal G}$ by erasing $v$ and $\mathrm{Star}(v)$ for all vertices of length $3$.
\end{itemize}

\item[(ii)] Assume $g\geq 2$ and let $M_0(D, N)_{st}$ denote the stable model of $X_0(D, N)$.

\begin{itemize}
\item If $p\ne 2, 3$ or $h_2=h_3=0$, $M_0(D, N)_{st}=M_0(D, N)$.

\item If $p=2$, $M_0(D, N)_{st}$ is the blow-down of all components $C_v$ of $M_0(D, N)$ with $\ell (v)=2$ or $3$.

\item If $p=3$, $M_0(D, N)_{st}$ is the blow-down of all components $C_v$ of $M_0(D, N)$ with $\ell(v)=2$ or $3$.

\end{itemize}

Its dual graph $\mathcal G_{st}$ is obtained from $\mathcal G$ by 

\begin{enumerate}

\item[-] If $p=2$, erasing $v$ and $\mathrm{Star}(v)$ for all $v$ with $\ell(v)=3$.

\item[-] If $p=2$ or $3$, replacing each chain $$v\stackrel{e}{-}v'\stackrel{e'}{-}v''$$ such that $\ell(v')=2$ or $3$, by $v\stackrel{(e'')}{-}v''$ with $\ell(e'')=\ell(e) +\ell(e')$. 
\end{enumerate}

\end{enumerate}

\end{proposition}

{\bf Proof.} $(i)\,$ By construction, $\tilde{M}_0(D, N)$ is regular. By Castelnuovo's criterion (cf.\,\cite[p.\,416-417]{Liu}), $\tilde{M}_0(D, N)$ is minimal over $\Z_p$ if and only if there exist no irreducible rational components $E$ in $\tilde{M}_0(D, N)_p$ which intersect the remaining components at a {\em single} point. This is equivalent to saying that there exists no vertex $v$ in $\tilde{\mathcal G}$ with $|\mathrm{Star}(v)|=1$. Since $|\mathrm{Star}(v)|=2$ for those vertices which were new created when constructing $\tilde{\mathcal G}$ from $\mathcal G$, we can directly work with $\mathcal G$. 

If $(D/p, N) = (2, 1)$ or $(3,1)$, the vertices $v$ and $v'$ are joined by $g+1\geq 2$ edges. Thus $\tilde{M}_0(D, N)$ is minimal.

Assume $(D/p, N) \ne (2, 1), (3,1)$. By (\ref{p+1}) and (\ref{star}), $|\mathrm{Star}(v)|=1$ for a vertex in $\mathcal G$ exactly when $p=2$ and $\ell(v)=3$. The minimal regular model is then obtained by blowing-down the corresponding components $C_v$.

$(ii)\,$ By definition, $M_0(D, N)$ is stable if $|\mathrm{Star}(v)|\geq 3$ for all vertices $v$ of $\mathcal G$. By (\ref{p+1}) and (\ref{star}), $|\mathrm{Star}(v)|<3$ exactly when $p=2$ or $3$ and $\ell (v)=2$ or $3$. When this holds, the stable model is achieved by blowing-down all these irreducible components. $\Box $

The incidence matrix of $\mathcal G$ can be recovered (and explicitly computed) from the theory of Brandt modules and matrices. Namely, 
let $M:=M_{\tcO ; p}\in \M_h(\Z )$ denote the Brandt matrix attached to $\tcO $ and the prime number $p$. Let $I_i$, $i=1, ..., h$, a set of representatives of left ideals of $\tcO $ up to principal ideals. By definition, $M(i, j)$ is the number of integral ideals of reduced norm $p$ which are equivalent on the right to $I_i^{-1}\cdot I_j$. 

Recall that there exist no edges joining vertices $v_i$, $v_j$ (and the same holds for $v_i'$, $v_j'$). The number of edges $e$ joining two given vertices $v_i$ and $v_j'$ (which equals that of edges joining $v_j$ and $v_i'$) can be computed by means of (\ref{star}) and the formula

\begin{equation}\label{nedges}
M(i, j) = \sum _{v_i\stackrel{e}{\ra} v_j'} \frac{\ell(v_i)}{\ell (e)}.
\end{equation}

In particular it always holds that $M(i, j)/\ell(v_i) = M(j, i)/\ell(v_j)$. This completely determines $\mathcal G$. Note that (\ref{p+1}) 
implies that $\sum_{j=1}^h M(i, j) = p+1$ for each row $i=1, ..., h$.

\subsection{Proof of Theorem \ref{CD}}\label{pcd}

\begin{definition}
An automorphism $\om \in \Aut (\mathcal G)$ is {\em admissible} if there exists no vertex $v$ in $\mathcal G $ fixed by $\om $ such that $\mathrm{Star}(v)$ contains three different edges $e_1, e_2, e_3$ also fixed by $\om $. A subgroup $A\subseteq \Aut (\mathcal G)$ is {\em admissible } if any $\om \in A$, $\om \ne Id $, is admissible.
\end{definition}

\begin{proposition}\label{admis}
Assume $g(\,X_0(D, N)\,) \geq 2$. Then there exists a monomorphism $\varrho :$ $\Aut (X_0(D, N) ) \hookrightarrow \Aut (\mathcal G)$ which embeds $\Aut (X_0(D, N))$ into an admissible subgroup of $\Aut (\mathcal G)$.
\end{proposition}

{\bf Proof.} By \cite[I.12]{DeMu}, \cite[Ch.\,IX]{Liu}, there exists a natural monomorphism $$\Aut (X_0(D, N) ) \hookrightarrow \Aut (M_0(D, N)_{st}\times \F_p).$$ 

Since $M_0(D, N)$ is the blow-up of $M_0(D, N)_{st}$ over $\Z_p$ at finitely many points, there is a birational morphism $M_0(D, N)\ra M_0(D, N)_{st}$ which induces an embedding $\Aut (M_0(D, N)_{st }\times \F_p) \hookrightarrow \Aut (M_0(D, N)\times \F_p)$.  

By considering the action on the irreducible components and singular points of the special fiber at $p$, any $\om \in \Aut (X_0(D, N))$ induces through the above embeddings an automorphism $\varrho (\om )$ of $\mathcal G$ as a graph with lengths. Let $\varrho :$ $\Aut (X_0(D, N) ) \hookrightarrow \Aut (\mathcal G)$ the resulting map. By construction, it is clearly a group homomorphism.

Assume that $\varrho (\om )=Id$. Then, the action of $\om $ on $M_0(D, N)_{st }\times \F_p$ would fix all its irreducible components and intersection points. Since a non-trivial automorphism of the projective line has at most two fixed points, we conclude that $\om = Id$. Hence $\varrho $ is a monomorphism and $\varrho (\om )$ is admissible for any $\om \ne Id$. $\Box $

As we mentioned above, the Atkin-Lehner involution $\omega_p$ acts on $\mathcal G$ as $\omega_p(v_i)=v_i'$. 

Proposition above provides an explicit method for proving in many instances that $s=r$. Indeed, the graph $\mathcal G$ can be computed by means of 
David Kohel's {\em Brandt modules} package implemented in MAGMA. Namely, for a given $p\mid D$, Brandt's matrix $M:=M_{\tcO ; p}$ can be computed by MAGMA V:=BrandtModule$(D/p,N)$; M:=HeckeOperator(V,p);. For a prime $q\mid D/p$, the action of the Atkin-Lehner involution $\om_q$ on $V$ is obtained by MAGMA $W_q$:=HeckeOperator(V,q);. The action of $\om_q$ on the set $E$ of edges of $\mathcal G$ is obtained by MAGMA E:=BrandtModule$(D/p,M p)$;\- $W'_q$:=HeckeOperator(E,q);.  

\vspace{0.2cm}

{\bf Proof of Theorem \ref{CD}.} Throughout we may assume that $(D,N)\ne (2 p, 1), (3 p,1)$ where $p$ denotes a prime number, as these cases are covered by Theorem \ref{all} (v).

$(i)$ Let $m=2$. Fix a prime $q\mid D$. Choose $q$ to be the (single) prime such that $(\frac{-4}{q}) = 1$ in case it exists. 
Let $\mathcal G$ the dual graph of $M_0(D, N)_q$. In $\mathrm{Ver}(\mathcal G)$ there exist $2\cdot h_2 = 2^{r-1}$ vertices of length $2$ if $2\nmid D N$ (resp.\,$2^{r-2}$ if $2\mid D N$). Since $q\ne 2$ and the maximal $2$-elementary subgroup of $\Aut(\mathbb{P}^1_{\F_{q^2}}) = \PGL_2(\F_{q^2})$ is isomorphic to $(\Z/2\Z )^2$ (cf.\,e.\,g.\,\cite[Lemma 1.3]{KR}), it follows that the stabilizer of any vertex has order at most $4$. 

Hence $s\leq r+1$ if $2\nmid D N$; $s=r$ if $2\mid D N$. The proof for $m=3$ follows along the same lines and one again obtains that if $(\frac{-3}{p}) \ne \varepsilon_p$ for all $p\mid D N$ except at most for one prime divisor of $D$ then $s\leq r+1$ if $3\nmid D N$; $s=r$ if $3\mid D N$. 

$(ii)$ Let $h=h(D/p, N)$. There is a well-defined action of $\Aut(X_0(D, N))/\langle \om_p\rangle $ on the subset $V=\{ v_1, ..., v_h\}$ of vertices of the dual graph $\mathcal G$ of $M_0(D, N)_p$. As above, the stabilizer of each of these vertices has order at most $4$ and the statement follows. 

$(iii)$ As we mentioned, there is a canonical embedding of $\Aut (X_0(D, N))$ into $\Aut (M_0(D, N)_{\ell })$. Since $\mathrm{ord }_2(0)=\infty $, we can assume that $M_0(D, N)_{\ell }(\F_{\ell }) \ne \emptyset $. Our claim now follows immediately from Lemma \ref{Og} applied to any point $P\in M_0(D, N)_p(\F_{\ell })$. $\Box $

\subsection{A numerical result}\label{nr}

Let us illustrate how our results above serve to prove that there exist no exceptional automorphisms in many Shimura curves.

\vspace{0.2cm}
{\bf Example.} Let $D=5\cdot 41$. We have $g(X_D)=13$. None of the items of Theorems \ref{CM} or \ref{CD} apply to show that all automorphisms of $X_{D}$ are modular. Since $h(41,1)=4$, the dual graph $\mathcal G$ of the special fiber of $X_{D}$ at $p=5$ consists of eight vertices $\{ v_1, ..., v_4, v_1', ..., v_4' \}$ of which $v_1$ and $v_1'$ have length $3$ while the remaining ones have length $1$. Moreover, all edges have length $1$. One computes that 
$$M=\begin{pmatrix}
0&3&0&3\\1&0&3&2\\0&3&2&1\\1&2&1&2
\end{pmatrix} \quad \text{ and } \quad W_{41} = Id_V.$$

In fact, the action of any involution $\om \in \Aut(\mathcal G)/\langle \om_5\rangle $ fixes each vertex $v\in V$, because $\mathrm{Star}\,(v_i)$ are pairwise different for $i=1, ..., 4$. Looking at the action of $\om $ on the rational component $C_{v_1}$, the two points $P, P'$ corresponding to the edges joining $v_1$ with $v_2'$ and $v_4'$ are necessarily fixed points of $\tilde w$. Since $\Aut(\tilde C_{v_1}, P, P') \simeq \Z/2\Z $, it follows that $\Aut(X_D) = \Aut (\mathcal G) = \langle \omega_5, \omega_{41}\rangle $.

\begin{proposition}\label{2000}
For $D\leq 1500$, the only automorphisms of $X_D$ are the Atkin-Lehner involutions, provided $g(X_D)\geq 2$ and $D\ne 493$, $583$, $667$, $697$, $943$.
\end{proposition}

{\bf Proof.} If the number of prime factors of $D$ is $r\geq 4$, Theorem \ref{CM} (i) and Corollary \ref{corol} show that $s=r$. 

Assume now that $D=p\cdot q$ is the product of two primes. By applying Theorems \ref{all} (v), \ref{CM} (i) and \ref{CD} (iii) for the three-hundred first primes $\ell \ne p, q$, we prove the statement for all such $D$ except for a list $L$ of $44$ values of $D$; we do not reproduce the list here for the sake of brevity. 

For the values of $D$ in $L$, we apply Proposition \ref{admis} as in the above example. This way we are able to prove that $s=2$ for all $D$ except for $D=85$, $145$, $493$, $583$, $667$, $697$, $943$. Let us illustrate what is going on with some examples. 

When $D=697$, let $\mathcal G$ denote the dual graph of the reduction mod $p=17$. With the notation above, it turns out that $V=\{ v_1, v_2, v_3, v_4\}$, with $\ell (v_1)=3$, $\ell (v_i)=1$ for $i=2, 3, 4$ and one computes that the permutation $\om $ of the vertices $v_2$ and $v_3$ is an admissible involution with commutes with the Atkin-Lehner involutions. Hence there exists an admissible subgroup of $\Aut (\mathcal G)$ which is isomorphic to $(\Z /2\Z )^3$ and we can not prove that $s=2$. 

When $D=1057$, let $\mathcal G$ denote the dual graph of the reduction mod $p=7$. One computes that $|V|=13$. Note that any automorphism of $X_{1057 }$ must induce a permutation $\om \in S_{13}$ of the vertices in $V$ which commutes with the matrixes $M$ and $W_{151}$. When computing $M$ and $W_{151}$, one shows that there exist exactly sixteen such permutations $\om_k$, $k=1, ..., 16$. But for each of them, it turns out that either $\om_k$ or $\om_k \cdot W_{151 }$ is {\em not} admissible. Thus $s=2$.  

It remain to prove our proposition for $D=85$, $145$. For $D=145$, it was already shown in the proof of \cite[Thm.\,7]{Rot:2002} that $\Aut(X_{145})\simeq \langle \om_5, \om_{29}\rangle $. One shows the same similarly for $D=85$: As obtained from William Stein's data basis, there exist isogenies $\mathrm{Jac}(X_{85}) \sim \mathrm{Jac}(X_0(85))^{new } \sim E\times S_1\times S_2$ defined over $\Q $, where $E$ is a (modular) elliptic curve and $S_1$, $S_2$ are modular abelian surfaces over $\Q $. 

The modularity implies that $\End (S_i)$ are orders in a real quadratic field. Hence, the only automorphisms of finite order of $E$, $S_1$ or $S_2$ are $\pm \mathrm{Id}$. Composing with this isogeny, we obtain a monomorphism $\Aut (X_{85}) \hookrightarrow A:=\{ \pm \mathrm{Id}_E \} \times \{ \pm \mathrm{Id}_{S_1}\} \times \{ \pm \mathrm{Id}_{S_2} \} $.  

By \cite{Ogg1}, $X_{85}$ is not hyperelliptic and this is saying that $(-\mathrm{Id}_E ,-\mathrm{Id}_{S_1}, -\mathrm{Id}_{S_2})$ does not lie in $\Aut(X_{85})$. Thus its index in $A$ is at least $2$ and we conclude that $s=2$. 

These ideas appear to be insufficient to prove the same result for $D= 493$, $583$, $667$, $697$ or $943$. For these $D$, we can just claim that $s\leq 3$, by Corollary \ref{Shimurabound}. $\Box $

\vspace{0.2cm}

{\bf Acknowledgement.} The authors would like to thank the referee for his useful remarks.

\end{document}